# Stein's method and Poisson process approximation for a class of Wasserstein metrics

DOMINIC SCHUHMACHER

*School of Mathematics and Statistics, University of Western Australia, 35 Stirling Highway, Crawley WA 6009, Australia. E-mail: [dominic@maths.uwa.edu.au](dominic@maths.uwa.edu.au)*

Based on Stein's method, we derive upper bounds for Poisson process approximation in the $L_1$-Wasserstein metric $d_2^{(p)}$, which is based on a slightly adapted $L_p$-Wasserstein metric between point measures. For the case $p = 1$, this construction yields the metric $d_2$ introduced in [Barbour and Brown *Stochastic Process. Appl.* **43** (1992) 9–31], for which Poisson process approximation is well studied in the literature. We demonstrate the usefulness of the extension to general $p$ by showing that $d_2^{(p)}$-bounds control differences between expectations of certain $p$th order average statistics of point processes. To illustrate the bounds obtained for Poisson process approximation, we consider the structure of 2-runs and the hard core model as concrete examples.

## 1. Introduction

Stein's method is a very powerful and flexible tool for deriving upper bounds for distances between probability distributions. Since its first publication in Stein (1972), where it was limited to normal approximation, the method has been extensively studied and adapted to a wide range of different distributions; see Barbour and Chen (2005) for a comprehensive overview. In Barbour and Brown (1992) (see also Barbour, Holst and Janson (1992) for discrete state spaces and the earlier results in Arratia, Goldstein and Gordon (1989) and Barbour (1988)) Poisson process approximation by Stein's method was developed both in the total variation metric and in a particular Wasserstein metric, denoted by $d_2$, that has proved to be more suitable for the problem of point process approximation. In Brown and Xia (2001) (after an earlier more complicated version in Brown, Weinberg and Xia (2000)) a partial improvement of the $d_2$-approximation was offered that was able to remove in many cases a rather annoying logarithmic factor from the upper bound. For a fine overview of Stein's method for Poisson process approximation see Xia (2005).

In the present paper we use Stein's method to give upper bounds for Poisson process approximation in a generalized $d_2$-metric, which we denote by $d_2^{(p)}$, where $p \in [1, \infty]$ and $d_2^{(1)} = d_2$. This generalization enables us to draw wider conclusions from the resulting estimates. In particular, we have that any upper bound obtained for a $d_2^{(p)}$-distance controls also the difference between the expectations of statistics that are based on the $p$th order average of certain distance features within the point processes, whereas often the same is true only for the standard (first-order) average in the case of $d_2$-bounds. The price to be paid for this improvement is that the upper bounds we obtain are in general somewhat worse. However, for $p < \infty$ they are still better than the corresponding total variation estimate, and they do not contain the infamous logarithmic factor that usually appears in the estimates for $p = 1$.

The paper is organized as follows. In Section 2 we give the definition of $d_2^{(p)}$ and discuss some of the elementary properties (Section 2.1). We furthermore present examples of the $p$th order average statistics mentioned above (Section 2.2). Section 3 contains our main result. After stating the general upper bound for Poisson process approximation in Section 3.1, we compute two examples in concrete situations (Section 3.2), before proving the bound in Section 3.3.

## 2. The Wasserstein metrics $d_2^{(p)}$

### 2.1. Notation and definitions

We always consider a compact metric space $(\mathcal{X}, d_0)$ with $d_0 \leq 1$ as the state space of our point processes and equip it with its Borel $\sigma$-algebra $\mathcal{B}$. Denote the space of all finite point measures on $\mathcal{X}$ by $\mathfrak{N}$ and equip it as usual with the vague topology and the $\sigma$-algebra $\mathcal{N}$ generated by this topology, which is the smallest $\sigma$-algebra that renders the point counts of measurable sets measurable (see Kallenberg (1986), Section 1.1, Lemma 4.1, and Section 15.7). Recall that a point process is just a random element of $\mathfrak{N}$.

We first define metrics $d_1^{(p)}$ on $\mathfrak{N}$ that are based on an $L_p$-Wasserstein construction. Denote the set of permutations of $\{1, 2, \ldots, n\}$ by $\Pi_n$. For any $\xi = \sum_{i=1}^{|\xi|} \delta_{x_i}$ and $\eta = \sum_{i=1}^{|\eta|} \delta_{y_i} \in \mathfrak{N}$, let

$$d_1^{(p)}(\xi, \eta) := \begin{cases} \min_{\pi \in \Pi_n} \left( \dfrac{1}{n} \sum_{i=1}^n d_0(x_i, y_{\pi(i)})^p \right)^{1/p}, & \text{if } |\xi| = |\eta| = n \geq 1, \\ 1, & \text{if } |\xi| \neq |\eta|, \\ 0, & \text{if } |\xi| = |\eta| = 0 \end{cases}$$

for $1 \leq p < \infty$, and let

$$d_1^{(\infty)}(\xi, \eta) := \begin{cases} \min_{\pi \in \Pi_n} \max_{1 \leq i \leq n} d_0(x_i, y_{\pi(i)}), & \text{if } |\xi| = |\eta| = n \geq 1, \\ 1, & \text{if } |\xi| \neq |\eta|, \\ 0, & \text{if } |\xi| = |\eta| = 0. \end{cases}$$

It is straightforward (in fact immediate, except for the triangle inequality, which can be proved by Minkowski's inequality) that $d_1^{(p)}$, $1 \leq p \leq \infty$, are metrics and bounded by 1. By applying Lyapunov's inequality, we obtain that

$$d_1^{(p)} \leq d_1^{(q)} \qquad \text{for } 1 \leq p \leq q \leq \infty, \tag{2.1}$$

and with the help of this result it can be seen that $d_1^{(p)}$ metrizes the vague topology for any $p$. Furthermore it can be shown that $(\mathfrak{N}, d_1^{(p)})$ is complete and separable (the latter follows directly from Result 15.7.7 in Kallenberg (1986)).

Next we define the metric $d_2^{(p)}$ on the space $\mathfrak{P}(\mathfrak{N})$ of probability measures on $\mathfrak{N}$, which is the usual $L_1$-Wasserstein metric with respect to $d_1^{(p)}$. Let $\mathcal{F}_2^{(p)} := \{f : \mathfrak{N} \to [0,1]; |f(\xi) - f(\eta)| \leq d_1^{(p)}(\xi, \eta) \text{ for all } \xi, \eta \in \mathfrak{N}\}$. Set then for any $P, Q \in \mathfrak{P}(\mathfrak{N})$

$$d_2^{(p)}(P, Q) := \sup_{f \in \mathcal{F}_2^{(p)}} \left| \int_{\mathfrak{N}} f \, \mathrm{d}P - \int_{\mathfrak{N}} f \, \mathrm{d}Q \right|.$$

Since this is exactly the Wasserstein construction (the fact that we restrict the functions in $\mathcal{F}_2^{(p)}$ to be $[0,1]$-valued has no influence on the supremum because the underlying $d_1^{(p)}$-metric is bounded by 1), it is clear that $d_2^{(p)}$, $1 \leq p \leq \infty$, are metrics, obviously bounded by 1, and that general results about Wasserstein metrics apply. One such result is the well-known Kantorovich–Rubinstein theorem, which in our situation states that

$$d_2^{(p)}(P, Q) = \min_{\substack{\Xi \sim P \\ \mathrm{H} \sim Q}} \mathbb{E} d_1^{(p)}(\Xi, \mathrm{H})$$

for any $P, Q \in \mathfrak{P}(\mathfrak{N})$, where we use notation of the form $Z \sim R$ to indicate that a random element $Z$ has distribution $R$. Furthermore it is clear by inequality (2.1) that

$$d_2^{(p)} \leq d_2^{(q)} \qquad \text{for } 1 \leq p \leq q \leq \infty, \tag{2.2}$$

and it follows, by the facts that $d_1^{(p)}$ metrizes the vague topology and that $d_2^{(p)}$ is also the *bounded* Wasserstein metric, that $d_2^{(p)}$ metrizes convergence in distribution of point processes (see Dudley (1989), Theorem 11.3.3).

To the author's knowledge, $d_2^{(p)}$ has not been considered before as a metric on $\mathfrak{P}(\mathfrak{N})$, except for $p = 1$ (as mentioned in the Introduction) and for $p = \infty$ (in Xia (1994) and Schuhmacher (2005a)).

### 2.2. Applications of distance estimates

By the definition of $d_2^{(p)}$, an upper bound of $d_2^{(p)}(\mathscr{L}(\Xi), \mathscr{L}(\mathrm{H}))$ controls also the difference $|\mathbb{E}f(\Xi) - \mathbb{E}f(\mathrm{H})|$ for any function $f \in \mathcal{F}_2^{(p)}$. It is thus of considerable interest in order to apply the upper bounds obtained in Theorem 3.A, to have a certain supply of "meaningful" $d_1^{(p)}$-Lipschitz continuous statistics of point patterns (where we do not worry too much about the Lipschitz constant as it will only appear as an additional factor in the upper bound). One way in which such statistics can then be used is to test if a given point pattern is a realization from among a certain class of point process distributions that are all known to lie within some $d_2^{(p)}$-distance $\varepsilon$ of a Poisson process distribution (e.g., according to our example in Section 3.2.2, the class of hard core processes with fixed intensity $\lambda$ and hard core radius $r$ below some level $\varrho > 0$). The fact that the test statistic lies in $\mathcal{F}_2^{(p)}$ enables us to control the size of the test in such a way that it lies only slightly below some required level $\alpha$ if $\varepsilon$ is small. A detailed application of this idea in the case $p = 1$ was presented in Schuhmacher (2005b), Section 3.2.

The examples of $d_1^{(p)}$-Lipschitz continuous statistics $f$ given below are all $p$th order averages of certain distance features within the point measure. In each case we tacitly set $f$ to zero where the stated definition does not apply (e.g. for $n < l$ in Proposition 2.A). The proofs of the propositions are given in the Appendix.

Our first example concerns $p$th order $U$-statistics with Lipschitz continuous kernels. Note that at least for $p = 1$ there is a plethora of results available about $U$-statistics that are based on a fixed number of i.i.d. points (which in the point process framework corresponds to a Poisson process conditioned on its total number of points). See Lee (1990) for more information. For $p = 1$, a class of functions similar to those in Proposition 2.A was proposed in Barbour, Holst and Janson (1992), Section 10.2.

**Proposition 2.A.** *Take $l \in \mathbb{N}$ and let $K: \mathbb{Z}_+ \times \mathcal{X}^l \to [0,1]$ be a function that is symmetric in the last $l$ arguments and satisfies*

$$|K(m; u_1, u_2, \ldots, u_l) - K(m; v_1, v_2, \ldots, v_l)| \leq \frac{1}{l} \sum_{i=1}^{l} d_0(u_i, v_i) \qquad (2.3)$$

*for all $m \in \mathbb{Z}_+$ and all $u_1, u_2, \ldots, u_l, v_1, v_2, \ldots, v_l \in \mathcal{X}$. Define $f: \mathfrak{N} \to [0,1]$ by*

$$f(\xi) := \overline{K^{(p)}}(\xi) := \left( \frac{1}{\binom{n}{l}} \sum_{1 \leq i_1 < i_2 < \cdots < i_l \leq n} K(n; x_{i_1}, x_{i_2}, \ldots, x_{i_l})^p \right)^{1/p} \qquad (2.4)$$

*for $\xi = \sum_{i=1}^{n} \delta_{x_i} \in \mathfrak{N}$ with $n \geq l$, and $1 \leq p < \infty$. Then $f \in \mathcal{F}_2^{(p)}$.*

Instead of (2.4), we may also consider the *centered* $p$th order average, which for the case $p = 2$ gives us the standard deviation of $(K(n; x_{i_1}, x_{i_2}, \ldots, x_{i_l}))_{1 \leq i_1 < i_2 < \cdots < i_l \leq n}$.

**Proposition 2.B.** *Let $K$ be as in Proposition 2.A and $\overline{K} := \overline{K^{(1)}}$. Define $f: \mathfrak{N} \to [0,1]$ by*

$$f(\xi) := \left( \frac{1}{\binom{n}{l}} \sum_{1 \leq i_1 < i_2 < \cdots < i_l \leq n} |K(n; x_{i_1}, x_{i_2}, \ldots, x_{i_l}) - \overline{K}(\xi)|^p \right)^{1/p} \qquad (2.5)$$

554	D. Schuhmacher

*for* $\xi = \sum_{i=1}^{n} \delta_{x_i} \in \mathfrak{N}$ *with* $n \geq l$, *and* $1 \leq p < \infty$. *Then* $f$ *is* $d_1^{(p)}$-*Lipschitz continuous with constant* 2.

One basic choice for the function $K$ in the above results is half the interpoint distance, that is, $K(m; u_1, u_2) = K_0(u_1, u_2) := \frac{1}{2} d_0(u_1, u_2)$ for all $m \in \mathbb{N}$ and $u_1, u_2 \in \mathcal{X}$. By the triangle inequality for $d_0$ it is immediately seen that inequality (2.3) holds. This choice allows several extensions to functions $K$ that are based on more than two points. One is half the average interpoint distance in groups of size $l \geq 2$, that is,

$$K_1(u_1, \ldots, u_l) := \frac{1}{2} \frac{1}{\binom{l}{2}} \sum_{1 \leq i < j \leq l} d_0(u_i, u_j).$$

Note that this function is only of interest for $p > 1$, since for $p = 1$ and any $l \geq 2$ we just obtain the same values $f(\xi)$ as under $K_0$ for any $\xi$ that has at least $l$ points. Let $\mathcal{X} \subset \mathbb{R}^D$, where for the sake of simplicity we assume that $\mathrm{diam}(\mathcal{X}) := \max\{|x - y|; x, y \in \mathcal{X}\} \leq 1$, and set $d_0(x, y) := |x - y|$. Then two more extensions are given as $2/l$ times the radius of the minimal bounding ball and $l/(2(l-1))$ times the average distance to the geometrical centroid (center of gravity) in groups of size $l$; that is, for $l \geq 2$ and $u_1, \ldots, u_l \in \mathcal{X}$,

$$K_2(u_1, \ldots, u_l) := \frac{2}{l} \min\{r \geq 0; \exists x \in \mathbb{R}^D \text{ such that } u_1, \ldots, u_l \in \mathbb{B}(x, r)\},$$

where $\mathbb{B}(x, r)$ denotes the closed Euclidean ball with center at $x$ and radius $r$, and

$$K_3(u_1, \ldots, u_l) := \frac{l}{2(l-1)} \frac{1}{l} \sum_{i=1}^{l} d_0\left(u_i, \frac{1}{l} \sum_{i=1}^{l} u_i\right).$$

For all of these functions $K_t$, $t \in \{1, 2, 3\}$, inequality (2.3) is straightforward to check by showing that

$$|K_t(u, u_2, \ldots, u_l) - K_t(v, u_2, \ldots, u_l)| \leq \frac{1}{l} d_0(u, v)$$

for all $u, v, u_2, \ldots, u_l \in \mathcal{X}$ and using the symmetry of $K_t$. More examples, some of which also have corresponding extensions to groups of size $l$, can be found in Schuhmacher (2005a).

Another $d_1^{(p)}$-Lipschitz continuous function is the $p$th order average of the nearest neighbor distances in a finite point measure on $\mathbb{R}^D$, where $D \in \mathbb{N}$. This statistic gives important information about the amount of clustering in a point pattern.

**Proposition 2.C.** *Let* $\mathcal{X} \subset \mathbb{R}^D$ *and* $d_0(x, y) := |x - y| \wedge 1$ *for all* $x, y \in \mathcal{X}$. *Define the function* $f : \mathfrak{N} \to [0, 1]$ *by*

$$f(\xi) := \left( \frac{1}{n} \sum_{i=1}^{n} \min_{\substack{j \in \{1, \ldots, n\} \\ j \neq i}} d_0(x_i, x_j)^p \right)^{1/p}$$

for $\xi = \sum_{i=1}^{n} \delta_{x_i} \in \mathfrak{N}$ *with* $n \geq 2$, *and* $1 \leq p < \infty$. *Then* $f$ *is* $d_1^{(p)}$-*Lipschitz continuous with constant* $\tau_D + 1$ *for* $p = 1$ *and* $2(2\tau_D + 1)^{1/p}$ *for general* $p$, *where* $\tau_D$ *denotes the kissing number in* $D$ *dimensions (i.e., the maximal number of unit balls that can touch a unit ball in* $(\mathbb{R}^D, |\cdot|)$ *without producing any overlaps of the interiors; see Conway and Sloane (1999), Section* 1.2, *for details).*

## 3. Distance bounds

In this subsection the main theorem is stated. We give an upper bound for $p \in [1, \infty]$ of the $d_2^{(p)}$-distance between the distribution of a general point process $\Xi$ and a Poisson process with the same expectation measure. The result is a generalization of Theorem 5.19 in Xia (2005) (the case $p = 1$), which in turn is ultimately based on Theorems 3.6 and 3.7 in Barbour and Brown (1992) (but incorporates among other things certain improvements made in Brown and Xia (1995a) and Chen and Xia (2004)).

### 3.1. Statement of the main theorem

We always consider a point process $\Xi$ on $\mathcal{X}$ that has finite expectation measure $\boldsymbol{\lambda}$. Let $\Xi_x$ be the Palm process of $\Xi$ given a point in $x$ (i.e. any point process that is distributed according to the Palm distribution of $\Xi$ given a point in $x$); see Kallenberg (1986), Chapter 10, for formal details or Xia (2005), Section 2.3.1, for a concise overview. Write $\lambda := |\boldsymbol{\lambda}|$ for the total mass of $\boldsymbol{\lambda}$, and denote by $\mathrm{Po}(\boldsymbol{\lambda})$ the Poisson process distribution with expectation measure $\boldsymbol{\lambda}$, and by $\mathrm{Po}(\lambda)$ the Poisson distribution with expectation $\lambda$.

Call a family $\{N_x\}_{x \in \mathcal{X}}$ of measurable subsets $N_x \subset \mathcal{X}$ a *neighborhood structure* if $x \in N_x$ and the mapping $[\mathfrak{N} \times \mathcal{X} \to \mathfrak{N}, (\xi, x) \mapsto \xi|_{N_x^c}]$ is $(\mathcal{N} \otimes \mathcal{B})$-$\mathcal{N}$-measurable. This is the case if $N(\mathcal{X}) := \{(x, y) \in \mathcal{X}^2; y \in N_x\}$ is $\mathcal{B}^2$-measurable (see Chen and Xia (2004), after Formula (2.4)). Note that $N_x$ does not have to be a neighborhood of $x$ in the topological sense.

If $\boldsymbol{\mu}$ is a finite measure on $\mathcal{X}$, then we say that the *density conditions are satisfied for* $\Xi$ (with respect to the reference measure $\boldsymbol{\mu}$) if $\Xi$ is a simple point process, and the Janossy densities $j_n : \mathcal{X}^n \to \mathbb{R}_+$ with respect to $\boldsymbol{\mu}^n$ exist for $n \geq 0$ and are hereditary (i.e., $j_n(x_1, \ldots, x_n) = 0$ implies $j_{n+1}(x_1, \ldots, x_n, x_{n+1}) = 0$ for all $x_1, \ldots, x_n, x_{n+1} \in \mathcal{X}$). In this case, it can be seen that a density $\phi : \mathcal{X} \to \mathbb{R}_+$ of the expectation measure $\boldsymbol{\lambda}$ with respect to $\boldsymbol{\mu}$ exists. Write furthermore $g(x; \xi)$ for the conditional density of having a point of $\Xi$ in $x$ given that $\Xi|_{N_x^c} = \xi$. See Xia (2005), Section 2.3.2, for details on Janossy densities and the definition of $g$, and see Schuhmacher (2008), Section 2.4 and Remark A.C, for the reason why hereditarity (or a similar property) is needed, as well as for an alternative approach using densities with respect to a Poisson process distribution rather than Janossy densities.

Define the metric $d_1'$ on $\mathfrak{N}$ by $d_1'(\xi, \eta) := (m - n) + \min_{\pi \in \Pi_m} \sum_{i=1}^{n} d_0(x_i, y_{\pi(i)})$ for $\xi = \sum_{i=1}^{n} \delta_{x_i}$ and $\eta = \sum_{i=1}^{m} \delta_{y_i}$ if $n \leq m$, and $d_1'(\xi, \eta) := d_1'(\eta, \xi)$ otherwise. Let $\kappa_0 :=$

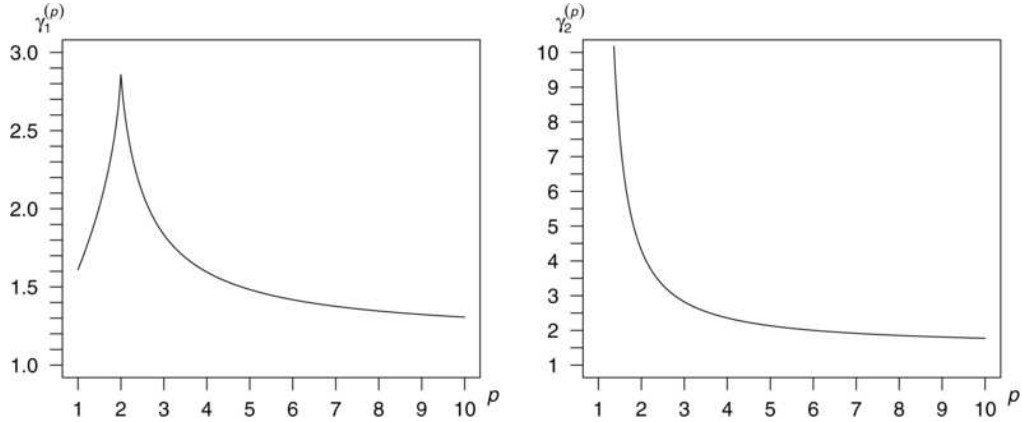

**Figure 1.** Graphs for $\gamma_1^{(p)}$ and $\gamma_2^{(p)}$. The limits for $p \to \infty$ are $1 + \frac{1}{4e(\kappa_0-1)}$ and $\frac{3}{2}$, respectively.

$4e/(1+4e-\sqrt{1+8e}) \approx 1.53$, $\kappa_1 := 2 - \kappa_0$,

$$\gamma_1^{(p)} := \begin{cases} \sqrt{\dfrac{2}{e}} + 2(\kappa_0 e^{\kappa_0})^{-1/2} \leq 1.61, & \text{if } p = 1, \\[2mm] \sqrt{\dfrac{2}{e}} + 2\left(\dfrac{2-p}{p\kappa_0^{1-1/p} - \kappa_1(p-1)}\right)^{(2-p)/(2(p-1))}, & \text{if } 1 < p \leq 2, \\[2mm] \dfrac{p}{p-1} + \dfrac{1}{\sqrt{2e}}\dfrac{p}{p-1}\left(\dfrac{p-2}{\sqrt{2e}(p\kappa_0^{1-1/p} - \kappa_1(p-1))}\right)^{(p-2)/p}, & \text{if } 2 < p < \infty, \end{cases}$$

and

$$\gamma_2^{(p)} := \frac{(1+2^{1/p}+(2/3)^{1/p})p^2}{(p-1)(2p-1)} \quad \text{for } 1 < p < \infty.$$

To get an impression of the behavior of $\gamma_1^{(p)}$ and $\gamma_2^{(p)}$ as functions of $p$, see Figure 1.

**Theorem 3.A.** *For any point process $\Xi$ on $\mathcal{X}$ with expectation measure $\boldsymbol{\lambda}$ and any neighborhood structure $(N_x)_{x \in \mathcal{X}}$, we have*

$$d_2^{(p)}(\mathscr{L}(\Xi), \mathrm{Po}(\boldsymbol{\lambda}))$$
$$\leq c_2^{(p)}(\lambda)\left(\int_{\mathcal{X}} \boldsymbol{\lambda}(N_x)\boldsymbol{\lambda}(\mathrm{d}x) + \mathbb{E}\int_{\mathcal{X}}(\Xi(N_x)-1)\Xi(\mathrm{d}x)\right) + \min(\varepsilon_1, \varepsilon_2),$$

*where*

$$\varepsilon_1 = c_1^{(p)}(\lambda)\mathbb{E}\int_{\mathcal{X}}|g(x;\Xi|_{N_x^c}) - \phi(x)|\boldsymbol{\mu}(\mathrm{d}x),$$

*which is valid if the density conditions are satisfied for $\Xi$ with respect to $\boldsymbol{\mu}$, and*

$$\varepsilon_2 = c_2^{(p)}(\lambda)\mathbb{E}\int_{\mathcal{X}} d_1'(\Xi|_{N_x^c}, \Xi_x|_{N_x^c})\boldsymbol{\lambda}(\mathrm{d}x).$$

*The factors $c_1^{(p)}(\lambda)$ and $c_2^{(p)}(\lambda)$ are given by*

$$c_1^{(p)}(\lambda) = \begin{cases} \min(1, \gamma_1^{(p)}\lambda^{-1/\max(2,p)}), & \text{if } 1 \leq p < \infty, \\ 1, & \text{if } p = \infty, \end{cases}$$

*and*

$$c_2^{(p)}(\lambda) = \begin{cases} \min(1, \frac{11}{6}[1 + 2\log^+(6\lambda/11)]\lambda^{-1}), & \text{if } p = 1, \\ \min(1, \gamma_2^{(p)}\lambda^{-1/p}), & \text{if } 1 < p < \infty, \\ 1, & \text{if } p = \infty. \end{cases}$$

***Remark 3.B.*** Note that $\gamma_2^{(p)} \to \infty$ for $p \to 1$, which is consistent with the fact that $c_2^{(1)}(\lambda)$ is not of the form "constant times $\lambda^{-1}$" but contains an extra factor of order $\log(\lambda)$. The presence of this factor in the upper bound of the $d_2^{(1)}$-distance has caused much discussion over the years, especially since no such factor is present in the corresponding upper bound of the total variation distance between the distributions of the total numbers of points (see Barbour and Brown (1992), Theorem 3.10).

It was shown in Brown and Xia (1995b) that with the current proof technique this factor cannot be omitted in a general setting (more precisely, that the estimate in Proposition 3.H(ii) is of the correct order if $p = 1$). In Brown et al. (2000) and Brown and Xia (2001) non-uniform bounds of the term $\Delta^2 h$ in Proposition 3.H(ii) were given, with the help of which the authors were able to dispose of the logarithmic factor in many important special cases. However, there is currently no general result available that can do without the logarithm. Very recently, Röllin (2008) gave an example of a point process $\Xi$, for which the (exact) order of $d_2^{(1)}(\mathscr{L}(\Xi), \mathrm{Po}(\boldsymbol{\lambda}))$ for $\lambda \to \infty$ contains an extra factor $\log(\lambda)$ as compared to the order of $d_{TV}(\mathscr{L}(|\Xi|), \mathrm{Po}(\lambda))$. This example makes the logarithmic term in $c_2^{(1)}(\lambda)$ appear rather natural.

## 3.2. Examples

In order to illustrate how the bounds given in Theorem 3.A can be used in concrete situations, we present two quick examples.

### 3.2.1. Process of 2-runs

This application has been considered for $p = 1$ in Section 6.2 of Xia (2005). The corresponding arguments remain largely the same.

Let $\mathcal{X} = [0, 1]$, $d_0 \leq 1$ an arbitrary metric on $\mathcal{X}$, and choose $0 < z_1 < z_2 < \cdots < z_n = 1$. Consider i.i.d. indicator random variables $I_1, I_2, \ldots, I_n$ with expectation $p$. In order to avoid edge effects, we interpret the indices $1, 2, \ldots, n$ as the elements of the quotient ring $\mathbb{Z}_n := \mathbb{Z}/n\mathbb{Z}$ (so that $n+1 = 1$ and $1-1 = n$). Define indicators $J_i := I_i I_{i+1}$ for $i \in \mathbb{Z}_n$. Then $\Xi := \sum_{i \in \mathbb{Z}_n} J_i \delta_{z_i}$ is a point process on $\mathcal{X}$ with expectation measure $\boldsymbol{\lambda} = \sum_{i \in \mathbb{Z}_n} p^2 \delta_{z_i}$, which describes the starting points of 2-runs in the process $\sum_{i \in \mathbb{Z}_n} I_i \delta_{z_i}$.

Applying Theorem 3.A is straightforward. Setting $N_{z_i} := \{z_i\}$, we can immediately see that

$$\int_{\mathcal{X}} \boldsymbol{\lambda}(N_x) \boldsymbol{\lambda}(\mathrm{d}x) = np^4 \quad \text{and} \quad \mathbb{E}\int_{\mathcal{X}} (\Xi(N_x) - 1)\Xi(\mathrm{d}x) = 0.$$

We give an upper bound for the term $\varepsilon_2$. As a concrete Palm process we may choose

$$\Xi_{z_i} = \delta_{z_i} + I_{i-1}\delta_{z_{i-1}} + I_{i+2}\delta_{z_{i+1}} + \sum_{j \in \mathbb{Z}_n \setminus \{i-1, i, i+1\}} J_j \delta_{z_j}.$$

For bounding $d_1'(\Xi|_{N_{z_i}^c}, \Xi_{z_i}|_{N_{z_i}^c})$, pair each point of $\Xi|_{N_{z_i}^c}$ with the corresponding point of $\Xi_{z_i}|_{N_{z_i}^c}$ at the same position, which gives a perfect match except at $z_{i-1}$ and $z_{i+1}$, where it can happen that $\Xi_{z_i}|_{N_{z_i}^c}$ has a point, but $\Xi|_{N_{z_i}^c}$ has none. Thus

$$d_1'(\Xi|_{N_{z_i}^c}, \Xi_{z_i}|_{N_{z_i}^c}) \leq I_{i-1} - J_{i-1} + I_{i+1} - J_{i+1} = I_{i-1}(1 - I_i) + I_{i+1}(1 - I_{i+2}),$$

which implies that

$$\varepsilon_2 \leq c_2^{(p)}(\lambda) 2np^3(1-p).$$

Collecting the various estimates, we obtain the following result.

**Proposition 3.C.** *With the above assumptions we have*

$$d_2^{(p)}(\mathcal{L}(\Xi), \mathrm{Po}(\boldsymbol{\lambda})) \leq \begin{cases} \frac{11}{6}[1 + 2\log^+(6np^2/11)] \cdot p(2-p), & \text{if } p = 1, \\ \gamma_2^{(p)}(np^2)^{1-1/p} \cdot p(2-p), & \text{if } 1 < p < \infty, \\ np^2 \cdot p(2-p), & \text{if } p = \infty. \end{cases}$$

*Remark 3.D.* In Theorem 6.3 of Xia (2005) it is shown that the logarithmic factor for $p = 1$ can be disposed of at the cost of a higher constant and a considerably more complicated proof.

*Remark 3.E.* The maybe more obvious choice of $N_{z_i} := \{z_{i-1}, z_i, z_{i+1}\}$ for the proof of Proposition 3.C, which implies that $\varepsilon_1 = \varepsilon_2 = 0$ in Theorem 3.A, would in fact yield a somewhat worse bound, where the factor $p(2-p)$ is replaced by $p(2+3p)$.

#### 3.2.2. Hard core process

This application has been considered for $p = 1$ in Barbour and Brown (1992) (see after Theorems 2.4 and 3.6), with an important correction in Brown and Greig (1994). The arguments below are largely the same as in the latter article.

Let $\mathcal{X} = [0,1]^D$ and $d_0 \leq 1$ an arbitrary metric on $\mathcal{X}$. In order to avoid edge effects, we shall assume that the torus convention holds, which will become important below when we measure Euclidean distances $|x - y|$. Let $\boldsymbol{\mu}$ be Lebesgue measure on $\mathcal{X}$, and consider a stationary hard core process $\Xi$ with expectation measure $\boldsymbol{\lambda} = \lambda\boldsymbol{\mu}$ for $\lambda > 0$ and with hard core radius $r > 0$ (note that $r$ cannot be above a certain threshold $r_0(\lambda) > 0$ that is determined by $\lambda$). Such a process may be specified by its Janossy densities with respect to $\boldsymbol{\mu}$, given by

$$j_n(x_1, \ldots, x_n) = c\beta^n \mathrm{I}[|x_i - x_j| > r \text{ for all } 1 \leq i < j \leq n],$$

where $c$ and $\beta$ are chosen in such a way that $\sum_{n=0}^\infty \int_{\mathcal{X}^n} (n!)^{-1} j_n(\mathbf{x}) \boldsymbol{\mu}^n(\mathrm{d}\mathbf{x}) = 1$ (correct normalization for $j_n$ to be Janossy densities) and $\sum_{n=0}^\infty \int_{\mathcal{X}^n} (n!)^{-1} j_{n+1}(x, \mathbf{y}) \boldsymbol{\mu}^n(\mathrm{d}\mathbf{y}) = \lambda$ for every $x \in \mathcal{X}$ (correct density of expectation measure, $\phi(x) \equiv \lambda$).

We can easily see that the density conditions are satisfied for $\Xi$, and we can thus apply Theorem 3.A and make use of the term $\varepsilon_1$. Setting $N_x := \{x\}$, it is immediately clear that the first two summands in the upper bound are zero. A short computation (see Brown and Greig (1994), Section 3) shows that $g(x; \xi) = \beta \mathrm{I}[\xi(\mathbb{B}(x,r)) = 0]$, where $\mathbb{B}(x,r)$ is the closed Euclidean ball with center at $x$ and radius $r$, and that $\mathbb{P}[\Xi(\mathbb{B}(x,r)) = 0] = \mathbb{P}[\Xi|_{N_x^c}(\mathbb{B}(x,r)) = 0] = \lambda/\beta$. By these two equations it can be easily seen that

$$\mathbb{E}\int_{\mathcal{X}} |g(x; \Xi|_{N_x^c}) - \phi(x)| \boldsymbol{\mu}(\mathrm{d}x) \leq 2\lambda\mathbb{E}(\Xi(\mathbb{B}(x,r))) = 2\lambda^2 \alpha_D r^D,$$

where $\alpha_D$ denotes the volume of $\mathbb{B}(0,1)$. Thus Theorem 3.A yields the following result.

**Proposition 3.F.** *With the above assumptions we have*

$$d_2^{(p)}(\mathscr{L}(\Xi), \mathrm{Po}(\boldsymbol{\lambda})) \leq \begin{cases} 2\gamma_1^{(p)} \alpha_D r^D \lambda^{2-1/\max(2,p)}, & \text{if } 1 \leq p < \infty, \\ 2\alpha_D r^D \lambda^2, & \text{if } p = \infty. \end{cases}$$

***Remark 3.G.*** Following the arguments in Section 4 of Brown and Greig (1994), it can be seen that the constant 2 in Proposition 3.F can be improved to 1.5 by choosing $N_x := \mathbb{B}(x, r/2)$, at the cost of an additional condition and a considerably more complicated proof.

### 3.3. Proof of Theorem 3.A

Stein's method for Poisson process approximation as originally developed in Barbour and Brown (1992) provides us with a general procedure for finding upper bounds for a distance term of the form $d(\mathscr{L}(\Xi), \mathrm{Po}(\boldsymbol{\lambda})) = \sup_{f \in \mathcal{F}} |\mathbb{E}f(\Xi) - \mathrm{Po}(\boldsymbol{\lambda})(f)|$ for some class $\mathcal{F}$ of measurable functions $f : \mathfrak{N} \to \mathbb{R}$.

The rough idea of this procedure is as follows. First, set up the so-called Stein equation as

$$f(\xi) - \mathrm{Po}(\boldsymbol{\lambda})(f) = \mathscr{A}h(\xi) \qquad \text{for } \xi \in \mathfrak{N}, \tag{3.1}$$

where $\mathscr{A}$ is given by

$$\mathscr{A}\tilde{h}(\xi) = \int_{\mathcal{X}} [\tilde{h}(\xi + \delta_x) - \tilde{h}(\xi)]\boldsymbol{\lambda}(\mathrm{d}x) + \int_{\mathcal{X}} [\tilde{h}(\xi - \delta_x) - \tilde{h}(\xi)]\xi(\mathrm{d}x)$$

for suitable functions $\tilde{h}: \mathfrak{N} \to \mathbb{R}$ and for $\xi \in \mathfrak{N}$. Thus $\mathscr{A}$ is the generator of the spatial immigration-death process with immigration measure $\boldsymbol{\lambda}$ and unit per capita death rate, for which $\mathrm{Po}(\boldsymbol{\lambda})$ plays the special role of being its stationary distribution (see Xia (2005), Section 3.2 for more information). Let $\mathbf{Z}_\xi$ be such an immigration-death process with starting configuration $\mathbf{Z}_\xi(0) = \xi \in \mathfrak{N}$. It can be shown that, if $f$ is bounded, the function $h = h_f : \mathfrak{N} \to \mathbb{R}$,

$$h(\xi) = h_f(\xi) := -\int_0^\infty [\mathbb{E}f(\mathbf{Z}_\xi(t)) - \mathrm{Po}(\boldsymbol{\lambda})(f)]\,\mathrm{d}t \qquad (3.2)$$

is well-defined and solves equation (3.1). Rather than bounding $|\mathbb{E}f(\Xi) - \mathrm{Po}(\boldsymbol{\lambda})(f)|$ directly, it is then the key idea of Stein's method to bound the equivalent term $|\mathbb{E}\mathscr{A}h(\Xi)|$, which in fact turns out to be a considerably easier task in many situations.

In Theorem 5.3 of Xia (2005), which is a (very slight) specialization of Theorem 2.3 in Chen and Xia (2004), this strategy is employed to give a very versatile but still somewhat raw upper bound, which incorporates the essence of several of the earlier results mentioned in the introduction. Note that we have interchanged $f$ and $h$ in our presentation, which results in notation that is more commonly used in the literature (see, e.g., Barbour, Holst and Janson (1992), Barbour and Brown (1992), or Brown and Xia (1995a)). A direct consequence of Theorem 5.3 is that, for any bounded measurable function $f: \mathfrak{N} \to \mathbb{R}_+$ and $h = h_f$ defined as in (3.2), we have

$$|\mathbb{E}f(\Xi) - \mathrm{Po}(\boldsymbol{\lambda})(f)| \qquad (3.3)$$
$$\leq \|\Delta^2 h\|_\infty \left( \int_{\mathcal{X}} \boldsymbol{\lambda}(N_x)\boldsymbol{\lambda}(\mathrm{d}x) + \mathbb{E}\int_{\mathcal{X}} (\Xi(N_x) - 1)\Xi(\mathrm{d}x) \right) + \min(\varepsilon_1(h), \varepsilon_2(h)),$$

where

$$\varepsilon_1(h) = \|\Delta h\|_\infty \mathbb{E} \int_{\mathcal{X}} |g(x; \Xi|_{N_x^c}) - \phi(x)|\boldsymbol{\mu}(\mathrm{d}x),$$

which is valid if the density conditions are satisfied for $\Xi$ with respect to $\boldsymbol{\mu}$, and

$$\varepsilon_2(h) = \mathbb{E}\int_{\mathcal{X}} |[h(\Xi|_{N_x^c} + \delta_x) - h(\Xi|_{N_x^c})] - [h(\Xi_x|_{N_x^c} + \delta_x) - h(\Xi_x|_{N_x^c})]|\boldsymbol{\lambda}(\mathrm{d}x).$$

Here, the supremum norms of the first and second differences of $h$ are defined as

$$\|\Delta h\|_\infty := \sup_{\xi \in \mathfrak{N}, x \in \mathcal{X}} |h(\xi + \delta_x) - h(\xi)|$$

and

$$\|\Delta^2 h\|_\infty := \sup_{\xi \in \mathfrak{N}; x,y \in \mathcal{X}} |h(\xi + \delta_x + \delta_y) - h(\xi + \delta_x) - h(\xi + \delta_y) + h(\xi)|.$$

Note that the above result does not make use of any particular metric $d$, since it does not restrict the choice of functions $f$ to a specific class $\mathcal{F}$. The refinement of the result by giving upper bounds on the various increments of $h = h_f$ according to special properties of $f$ is the crucial step in adapting Stein's method to any one particular metric and is typically quite complicated. This step for the metrics $d_2^{(p)}$ is made in Proposition 3.H below. Inequality (3.3) together with this proposition directly yields the statement of Theorem 3.A.

**Proposition 3.H.** *Let $p \in [1, \infty]$. If $f \in \mathcal{F}_2^{(p)}$, then*

  (i) $\|\Delta h\|_\infty \leq c_1^{(p)}(\lambda)$;
  (ii) $\|\Delta^2 h\|_\infty \leq c_2^{(p)}(\lambda)$;
  (iii) $|(h(\xi + \delta_x) - h(\xi)) - (h(\eta + \delta_x) - h(\eta))| \leq c_2^{(p)}(\lambda) d_1'(\xi, \eta)$ *for $\xi, \eta \in \mathfrak{N}$ and $x \in \mathcal{X}$.*

**Proof.** The proof builds on the ideas of the proofs of the corresponding results for the case $p = 1$; see Propositions 5.16 to 5.18 in Xia (2005). In particular, it makes use of the representation of the spatial immigration-death process $\mathbf{Z}_\xi$ as $\mathbf{Z}_\xi(t) \stackrel{\mathscr{D}}{=} \mathbf{D}_\xi(t) + \mathbf{Z}_0(t)$, where $\mathbf{D}_\xi$ is a spatial pure death process with unit per capita death rate and starting configuration $\xi$, $\mathbf{Z}_0$ is a spatial immigration-death process with the same parameters as $\mathbf{Z}_\xi$, but starting with 0-measure, and $\mathbf{D}_\xi$ and $\mathbf{Z}_0$ are independent (see Xia (2005), Proposition 3.5). Write $Z_{|\xi|}(t) := |\mathbf{Z}_\xi(t)|$, $Z_0(t) := |\mathbf{Z}_0(t)|$, and note that $Z_0(t) \sim \mathrm{Po}(\lambda_t)$, where $\lambda_t = \lambda(1 - \mathrm{e}^{-t})$.

*Statement* (i). Suppose that $1 < p < \infty$. Inequality (5.19) in Xia (2005) yields that

$$|h(\xi + \delta_x) - h(\xi)| \leq \int_0^\infty \mathrm{e}^{-t} \{1 \wedge [|\mathbb{E}(f(\mathbf{Z}_\xi(t) + \delta_x) - f(\mathbf{Z}_\xi(t) + \delta_U))| \\ + |\mathbb{E}(f(\mathbf{Z}_\xi(t) + \delta_U) - f(\mathbf{Z}_\xi(t)))|]\} \, \mathrm{d}t \qquad (3.4)$$

for a random element $U \sim \boldsymbol{\lambda}/\lambda$ of $\mathcal{X}$ that is independent of everything else, where

$$\begin{aligned}
|\mathbb{E}(f(\mathbf{Z}_\xi(t) + \delta_x) - f(\mathbf{Z}_\xi(t) + \delta_U))| &\leq \mathbb{E}\left(\left(\frac{1}{|\mathbf{Z}_\xi(t)| + 1}\right)^{1/p}\right) \\
&\leq \left(\mathbb{E}\left(\frac{1}{Z_{|\xi|}(t) + 1}\right)\right)^{1/p} \\
&\leq \left(\mathbb{E}\left(\frac{1}{Z_0(t) + 1}\right)\right)^{1/p} \\
&= \left(\frac{1 - \mathrm{e}^{-\lambda_t}}{\lambda_t}\right)^{1/p}
\end{aligned} \qquad (3.5)$$

(see inequality (3.11) for details on the first estimate), and

$$|\mathbb{E}(f(\mathbf{Z}_\xi(t) + \delta_U) - f(\mathbf{Z}_\xi(t)))| \leq \frac{1}{\sqrt{2e\lambda_t}}$$

by inequality (5.23) in Xia (2005) (note that "=" in the last line should be "≤"). Hence,

$$
\begin{aligned}
|h(\xi + \delta_x) - h(\xi)| \\
&\leq \int_0^\infty e^{-t} \left\{ 1 \wedge \left[ \left( \frac{1 - e^{-\lambda_t}}{\lambda_t} \right)^{1/p} + \frac{1}{\sqrt{2e\lambda_t}} \right] \right\} dt \\
&= \int_0^1 \left\{ 1 \wedge \left[ \left( \frac{1 - e^{-\lambda s}}{\lambda s} \right)^{1/p} + \frac{1}{\sqrt{2e\lambda s}} \right] \right\} ds \qquad (3.6) \\
&\leq \frac{\kappa_0}{\lambda} + \int_{\kappa_0/\lambda}^1 \left[ \left( \frac{1}{\lambda s} \right)^{1/p} + \frac{1}{\sqrt{2e\lambda s}} \right] ds \\
&= \frac{\kappa_0}{\lambda} + \frac{p}{p-1} \left( \frac{1}{\lambda} \right)^{1/p} \left( 1 - \left( \frac{\kappa_0}{\lambda} \right)^{1-1/p} \right) + \sqrt{\frac{2}{e}} \frac{1}{\sqrt{\lambda}} \left( 1 - \sqrt{\frac{\kappa_0}{\lambda}} \right)
\end{aligned}
$$

for $\lambda \geq \kappa_0$, where $\kappa_0$ was defined such that it satisfies $\kappa_0^{-1} + (2e\kappa_0)^{-1/2} = 1$. Write $\kappa(p) := \frac{p}{p-1} \kappa_0^{1-1/p} - \kappa_1$, which can be easily seen to be strictly decreasing in $p$ with limit $2(\kappa_0 - 1) > 0$ for $p \to \infty$. For $1 < p \leq 2$, we factor out $\lambda^{-1/2}$, and maximize the left-over term

$$\sqrt{\frac{2}{e}} + \frac{p}{p-1} \left( \frac{1}{\lambda} \right)^{1/p - 1/2} - \kappa(p) \left( \frac{1}{\lambda} \right)^{1/2} \qquad (3.7)$$

in $\lambda$. For $1 < p < 2$, taking the first and second derivatives shows that a global maximum is attained at

$$\sqrt{\frac{2}{e}} + \frac{p}{p-1} \left( \frac{2-p}{(p-1)\kappa(p)} \right)^{(2-p)/(2(p-1))} - \kappa(p) \left( \frac{2-p}{(p-1)\kappa(p)} \right)^{p/(2(p-1))} = \gamma_1^{(p)}.$$

For $p = 2$, the term (3.7) is obviously strictly increasing in $\lambda$, so that letting $\lambda \to \infty$ yields that $\gamma_1^{(p)}$ maximizes this term also in the case $p = 2$. Thus, by inequality (3.6), $\|\Delta h\|_\infty \leq \gamma_1^{(p)} \lambda^{-1/2}$ for $1 < p \leq 2$.

For $p > 2$, we factor out $\lambda^{-1/p}$ in inequality (3.6), and maximize the left-over term

$$\frac{p}{p-1} + \sqrt{\frac{2}{e}} \left( \frac{1}{\lambda} \right)^{1/2 - 1/p} - \kappa(p) \left( \frac{1}{\lambda} \right)^{1 - 1/p}$$

in $\lambda$. Taking the first and second derivatives shows that a global maximum is attained at

$$\frac{p}{p-1} + \sqrt{\frac{2}{e} \left( \frac{1}{2e} \left( \frac{p-2}{(p-1)\kappa(p)} \right)^2 \right)^{1/2 - 1/p}} - \kappa(p) \left( \frac{1}{2e} \left( \frac{p-2}{(p-1)\kappa(p)} \right)^2 \right)^{1 - 1/p} = \gamma_1^{(p)}.$$

Thus, by inequality (3.6), $\|\Delta h\|_\infty \leq \gamma_1^{(p)} \lambda^{-1/p}$ for $p > 2$.

In total we have shown, for $1 < p < \infty$, that $\|\Delta h\|_\infty \leq \gamma_1^{(p)} \lambda^{-1/\max(2,p)}$ if $\lambda \geq \kappa_0$. By equation (3.4) we have $\|\Delta h\|_\infty \leq \int_0^\infty e^{-t} dt = 1$ for any $\lambda$. Statement (i) is then obtained because $\gamma_1^{(p)} \geq \kappa_0^{1/\max(2,p)} > \lambda^{1/\max(2,p)}$ if $\lambda < \kappa_0$, which follows for $p > 2$ simply by $\frac{p}{p-1} \geq \kappa_0^{1/p}$ and for $1 < p \leq 2$ by using the alternative expression via $x$ from (3.8), the inequality $(1+y)^r < \exp(ry)$ for $r, y > 0$, and that $(x+2)\kappa_0^{1/(x+2)} - \kappa_1 - x$ is maximal at $x = 0$.

What remains to be shown are the cases $p = 1$ and $p = \infty$. Since $\|\Delta h\|_\infty \leq 1$ holds always, the statement for $p = \infty$ is clear. For $p = 1$, we make use of the fact that $f \in \mathcal{F}_2^{(1)}$ implies $f \in \mathcal{F}_2^{(p)}$ and thus $|h(\xi + \delta_x) - h(\xi)| \leq c_1^{(p)}(\lambda)$ holds for every $p > 1$. Letting $p \to 1$ yields the required upper bound, where $\gamma_1^{(p)} \to \gamma_1^{(1)}$ follows by substituting $x := \frac{2-p}{p-1}$, so that

$$\left(\frac{2-p}{(p-1)\kappa(p)}\right)^{-(2-p)/(p-1)} = \left(1 + \frac{(x+2)\kappa_0^{1/(x+2)} - \kappa_1 - x}{x}\right)^x$$
$$\longrightarrow \exp(2 - \kappa_1 + \log(\kappa_0)) = \kappa_0 e^{\kappa_0} \quad \text{as } x \to \infty. \tag{3.8}$$

*Statement* (ii). Suppose that $1 < p < \infty$. As in the first part of the proof of Proposition 5.17 in Xia (2005), we obtain that

$$h(\xi + \delta_x + \delta_y) - h(\xi + \delta_x) - h(\xi + \delta_y) + h(\xi)$$
$$= -\int_0^\infty e^{-2t} \mathbb{E}[f(\mathbf{Z}_\xi(t) + \delta_x + \delta_y) - f(\mathbf{Z}_\xi(t) + \delta_x) \tag{3.9}$$
$$- f(\mathbf{Z}_\xi(t) + \delta_y) + f(\mathbf{Z}_\xi(t))] dt,$$

where there are numbers $b_k(t) \in [0, 1]$ for $k \geq -1$ such that

$$\mathbb{E}[f(\mathbf{Z}_\xi(t) + \delta_x + \delta_y) - f(\mathbf{Z}_\xi(t) + \delta_x) - f(\mathbf{Z}_\xi(t) + \delta_y) + f(\mathbf{Z}_\xi(t))]$$
$$\leq \mathbb{E}\left(\left(\frac{1}{Z_{|\xi|}(t) + 2}\right)^{1/p}\right) + \mathbb{E}\left(\left(\frac{\mathrm{I}[Z_{|\xi|}(t) \geq 1]}{Z_{|\xi|}(t)}\right)^{1/p}\right) \tag{3.10}$$
$$+ \sum_{k=-1}^\infty b_k(t)(\mathbb{P}[Z_0(t) = k-1] - 2\mathbb{P}[Z_0(t) = k] + \mathbb{P}[Z_0(t) = k+1]).$$

The only difference between (3.10) and the corresponding inequality on page 155 of Xia (2005) are the exponents $1/p$. They stem from a straightforward adaptation of inequalities (5.24) and (5.26) in Xia (2005) (note that "=" in the last line of (5.26) should be "$\leq$"), which is obtained by employing the estimate

$$|f(\eta + \delta_x) - f(\eta + \delta_y)| \leq d_1^{(p)}(\eta + \delta_x, \eta + \delta_y) \leq \left(\frac{1}{|\eta|+1}\right)^{1/p} d_0(x, y) \leq \left(\frac{1}{|\eta|+1}\right)^{1/p} \tag{3.11}$$

for $\eta \in \mathfrak{N}$. Continuing from equation (3.10), we have

$$\sum_{k=-1}^{\infty} b_k(t)(\mathbb{P}[Z_0(t) = k-1] - 2\mathbb{P}[Z_0(t) = k] + \mathbb{P}[Z_0(t) = k+1]) \leq \frac{1}{\lambda_t} \quad (3.12)$$

as shown on page 155 in Xia (2005), and

$$\mathbb{E}\left(\left(\frac{1}{Z_{|\xi|}(t) + 2}\right)^{1/p}\right) + \mathbb{E}\left(\left(\frac{\mathrm{I}[Z_{|\xi|}(t) \geq 1]}{Z_{|\xi|}(t)}\right)^{1/p}\right)$$
$$\leq \left(2^{1/p} + \left(\frac{2}{3}\right)^{1/p}\right) \mathbb{E}\left(\left(\frac{1}{Z_{|\xi|}(t) + 1}\right)^{1/p}\right) \quad (3.13)$$
$$\leq \frac{2^{1/p} + (2/3)^{1/p}}{\lambda_t^{1/p}},$$

where the first inequality is obtained because the sequence $((\frac{k+1}{k+2})^{1/p} + (\frac{k+1}{k})^{1/p})_{k \in \mathbb{N}}$ is seen to be decreasing, and the second inequality follows from (3.5).

In total, we combine (3.9), (3.10), (3.12), and (3.13), replacing $f$ by $(1-f) \in \mathcal{F}_2^{(p)}$ in (3.10) if necessary, to obtain

$$|h(\xi + \delta_x + \delta_y) - h(\xi + \delta_x) - h(\xi + \delta_y) + h(\xi)|$$
$$\leq \int_0^{\infty} \mathrm{e}^{-2t} \left\{2 \wedge \left[\frac{2^{1/p} + (2/3)^{1/p}}{\lambda_t^{1/p}} + \frac{1}{\lambda_t}\right]\right\} \mathrm{d}t$$
$$= \int_0^1 (1-s) \left\{2 \wedge \left[(2^{1/p} + (2/3)^{1/p})\left(\frac{1}{\lambda s}\right)^{1/p} + \frac{1}{\lambda s}\right]\right\} \mathrm{d}s \quad (3.14)$$
$$\leq 2\frac{\kappa_2(p)}{\lambda} - \left(\frac{\kappa_2(p)}{\lambda}\right)^2 + \int_{\kappa_2(p)/\lambda}^1 (1-s)\beta(p)\left(\frac{1}{\lambda s}\right)^{1/p} \mathrm{d}s$$
$$= 2\frac{\kappa_2(p)}{\lambda} - \left(\frac{\kappa_2(p)}{\lambda}\right)^2$$
$$+ \beta(p)\left(\frac{1}{\lambda}\right)^{1/p}\left[\frac{p^2}{(p-1)(2p-1)} - \frac{p}{p-1}\left(\frac{\kappa_2(p)}{\lambda}\right)^{1-1/p} + \frac{p}{2p-1}\left(\frac{\kappa_2(p)}{\lambda}\right)^{2-1/p}\right]$$

for $\lambda \geq \kappa_2(p)$, where $\kappa_2(p) := (\beta(p)/2)^p$ and $\beta(p) := (1 + 2^{1/p} + (2/3)^{1/p})$. We factor out $\lambda^{-1/p}$, and find a bound for the left-over term

$$\frac{\beta(p)p^2}{(p-1)(2p-1)} - \left(\frac{\beta(p)}{2}\right)^p \frac{2}{p-1}\left(\frac{1}{\lambda}\right)^{1-1/p} + \left(\frac{\beta(p)}{2}\right)^{2p} \frac{1}{2p-1}\left(\frac{1}{\lambda}\right)^{2-1/p}$$

on $\lambda \in [(\beta(p)/2)^p, \infty)$. From the first derivative we can see that this term is strictly increasing on the whole interval, so that the desired bound is obtained by letting $\lambda$ go to infinity. Hence $\|\Delta^2 h\|_\infty \leq \gamma_2^{(p)} \lambda^{-1/p}$ if $\lambda \geq (\beta(p)/2)^p$ and, by the first inequality in (3.14), $\|\Delta^2 h\|_\infty \leq \int_0^\infty 2e^{-2t}\,dt = 1$ for any $\lambda$. Noting that $\gamma_2^{(p)} \lambda^{-1/p} > 1$ for $\lambda < (\beta(p)/2)^p$, we obtain Statement (ii) for $1 < p < \infty$.

The case $p = 1$ was proved as Proposition 5.17 in Xia (2005). Since $\|\Delta^2 h\|_\infty \leq 1$ holds always, the case $p = \infty$ is obvious.

*Statement* (iii). Suppose that $1 < p < \infty$. We step by step adapt the proof of Proposition 5.18 in Xia (2005). Write $\xi = \sum_{i=1}^n \delta_{x_i}$ and $\eta = \sum_{i=1}^m \delta_{y_i}$, assuming without loss of generality that $n \leq m$ and that the points of $\xi$ and $\eta$ are numbered according to a $d_1'$-pairing, that is, such that $(m-n) + \sum_{i=1}^n d_0(x_i, y_i) = d_1'(\xi, \eta)$. Let $\eta_j := \sum_{i=1}^{n+j} \delta_{y_i}$ for $0 \leq j \leq m - n$. Then

$$|(h(\xi + \delta_x) - h(\xi)) - (h(\eta + \delta_x) - h(\eta))|$$
$$\leq |(h(\xi + \delta_x) - h(\xi)) - (h(\eta_0 + \delta_x) - h(\eta_0))| \qquad (3.15)$$
$$+ |(h(\eta_0 + \delta_x) - h(\eta_0)) - (h(\eta + \delta_x) - h(\eta))|,$$

where the second summand can be estimated as

$$|(h(\eta_0 + \delta_x) - h(\eta_0)) - (h(\eta + \delta_x) - h(\eta))|$$
$$\leq \sum_{j=1}^{m-n} |(h(\eta_j + \delta_x) - h(\eta_j)) - (h(\eta_{j-1} + \delta_x) - h(\eta_{j-1}))| \qquad (3.16)$$
$$\leq \|\Delta^2 h\|_\infty (m - n).$$

The first summand in (3.15) is zero if $n = 0$. For $n \geq 1$, write $\xi_j = \sum_{i=1}^{j-1} \delta_{x_i} + \sum_{i=j}^n \delta_{y_i}$ for $1 \leq j \leq n+1$, so that

$$|(h(\xi + \delta_x) - h(\xi)) - (h(\eta_0 + \delta_x) - h(\eta_0))|$$
$$\leq \sum_{j=1}^n |(h(\xi_{j+1} + \delta_x) - h(\xi_{j+1})) - (h(\xi_j + \delta_x) - h(\xi_j))|$$
$$\leq \sum_{j=1}^n \left( d_0(x_j, y_j) \int_0^\infty 2e^{-2t} \mathbb{E}\left( \left( \frac{1}{Z_{n-1}(t) + 1} \right)^{1/p} \right) dt \right) \qquad (3.17)$$
$$\leq \left( \int_0^\infty e^{-2t} \left\{ 2 \wedge \frac{2^{1/p} + (2/3)^{1/p}}{\lambda_t^{1/p}} \right\} dt \right) \sum_{j=1}^n d_0(x_j, y_j)$$
$$\leq c_2^{(p)}(\lambda) \sum_{j=1}^n d_0(x_j, y_j),$$

where the second estimate is obtained from the first inequality in the proof of Lemma 5.15 in Xia (2005) (adjusted by using (3.11) in the last line), the third estimate holds by inequality (3.5), and the last estimate follows from the proof of Statement (ii) (which shows that the second line of inequality (3.14) is bounded by $c_2^{(p)}(\lambda)$). The combining of (3.15), (3.16), and (3.17) yields Statement (iii) for $1 < p < \infty$.

The case $p = 1$ was proved as Propositon 5.18 in Xia (2005). The case $p = \infty$ follows by the same proof as above, but bounding the term in the second line of inequality (3.17) by $\sum_{j=1}^{n}(d_0(x_j, y_j) \int_0^\infty 2e^{-2t}\,dt) = \sum_{j=1}^{n} d_0(x_j, y_j)$, which is done by using $|f(\eta + \delta_x) - f(\eta + \delta_y)| \leq d_0(x, y)$ instead of (3.11) in the proof of Proposition 5.15 in Xia (2005). $\square$

# Appendix: Proofs of the Lipschitz continuities in Section 2.2

**Proof of Proposition 2.A.** It is obvious that $|f(\xi) - f(\eta)| \leq d_1^{(p)}(\xi, \eta)$ is satisfied for $\xi, \eta \in \mathfrak{N}$ with $|\xi| \neq |\eta|$ (since $\mathrm{im}(f) \subset [0, 1]$) or with $|\xi| = |\eta| < l$ (since in this case $f(\xi) = f(\eta) = 0$). Suppose then that $\xi = \sum_{i=1}^{n} \delta_{x_i}$ and $\eta = \sum_{i=1}^{n} \delta_{y_i}$, where $n \geq l$ and where the points of $\xi$ and $\eta$ are numbered according to a $d_1^{(p)}$-pairing, that is, such that $(\frac{1}{n}\sum_{i=1}^{n} d_0(x_i, y_i)^p)^{1/p} = d_1^{(p)}(\xi, \eta)$. Note that inequality (2.3) together with Lyapunov's inequality implies that

$$|K(m, u_1, \ldots, u_l) - K(m, v_1, \ldots, v_l)|^p \leq \left(\frac{1}{l}\sum_{i=1}^{l} d_0(u_i, v_i)\right)^p \leq \frac{1}{l}\sum_{i=1}^{l} d_0(u_i, v_i)^p.$$

Using the inverse triangle inequality for $\ell_p$-norms in the first line, we then obtain that

$$\begin{aligned}
|f(\xi) - f(\eta)|^p &\leq \frac{1}{\binom{n}{l}} \sum_{1 \leq i_1 < \cdots < i_l \leq n} |K(n; x_{i_1}, \ldots, x_{i_l}) - K(n; y_{i_1}, \ldots, y_{i_l})|^p \\
&\leq \frac{1}{l} \frac{1}{\binom{n}{l}} \sum_{1 \leq i_1 < \cdots < i_l \leq n} \sum_{r=1}^{l} d_0(x_{i_r}, y_{i_r})^p \\
&= \frac{1}{l} \frac{1}{\binom{n}{l}} \binom{n-1}{l-1} \sum_{i=1}^{n} d_0(x_i, y_i)^p \\
&= (d_1^{(p)}(\xi, \eta))^p. \qquad \square
\end{aligned} \tag{A.1}$$

**Proof of Proposition 2.B.** We show $|f(\xi) - f(\eta)| \leq 2d_1^{(p)}(\xi, \eta)$ in the non-obvious case. Let $\xi = \sum_{i=1}^{n} \delta_{x_i}$ and $\eta = \sum_{i=1}^{n} \delta_{y_i}$, where again $n \geq l$ and the points of $\xi$ and $\eta$ are numbered according to a $d_1^{(p)}$-pairing. Using the inverse triangle inequality for $\ell_p$-norms for the first and the usual triangle inequality for the second relation, we obtain that

$$|f(\xi) - f(\eta)|$$
$$\leq \left(\frac{1}{\binom{n}{l}} \sum_{1 \leq i_1 < \cdots < i_l \leq n} |(K(n; x_{i_1}, \ldots, x_{i_l}) - K(n; y_{i_1}, \ldots, y_{i_l})) - (\overline{K}(\xi) - \overline{K}(\eta))|^p\right)^{1/p}$$
$$\leq \left(\frac{1}{\binom{n}{l}} \sum_{1 \leq i_1 < \cdots < i_l \leq n} |K(n; x_{i_1}, \ldots, x_{i_l}) - K(n; y_{i_1}, \ldots, y_{i_l})|^p\right)^{1/p} + |\overline{K}(\xi) - \overline{K}(\eta)|$$
$$\leq 2d_1^{(p)}(\xi, \eta)$$

by inequality (A.1) (once for general $p$ and once for $p = 1$) and inequality (2.1). □

**Proof of Proposition 2.C.** Obviously, $|f(\xi) - f(\eta)| \leq d_1^{(p)}(\xi, \eta)$ if $|\xi| \neq |\eta|$ or $|\xi| = |\eta| < 2$. Suppose then that $\xi = \sum_{i=1}^n \delta_{x_i}$ and $\eta = \sum_{i=1}^n \delta_{y_i}$, where $n \geq 2$ and the points of $\xi$ and $\eta$ are numbered according to a $d_1^{(p)}$-pairing. Let $J(i)$ be the index of a nearest neighbor (with respect to $|\cdot|$ and hence $d_0$) of $x_i$ within the points of $\xi$ and $K(i)$ the index of a nearest neighbor of $y_i$ within the points of $\eta$. For $i$ fixed, we have

$$d_0(x_i, x_{J(i)}) \leq d_0(x_i, x_{K(i)}) \leq d_0(x_i, y_i) + d_0(y_i, y_{K(i)}) + d_0(y_{K(i)}, x_{K(i)}),$$

and

$$d_0(y_i, y_{K(i)}) \leq d_0(y_i, y_{J(i)}) \leq d_0(y_i, x_i) + d_0(x_i, x_{J(i)}) + d_0(x_{J(i)}, y_{J(i)}),$$

so that altogether

$$|d_0(x_i, x_{J(i)}) - d_0(y_i, y_{K(i)})| \leq d_0(x_i, y_i) + d_0(x_{L(i)}, y_{L(i)}),$$

where $L(i) := K(i)$ if $d_0(x_i, x_{J(i)}) \geq d_0(y_i, y_{K(i)})$ and $L(i) := J(i)$ otherwise. By the inverse triangle inequality for $\ell_p$-norms, we obtain now

$$|f(\xi) - f(\eta)|^p \leq \frac{1}{n} \sum_{i=1}^n |d_0(x_i, x_{J(i)}) - d_0(y_i, y_{K(i)})|^p$$
$$\leq \frac{1}{n} \sum_{i=1}^n (d_0(x_i, y_i) + d_0(x_{L(i)}, y_{L(i)}))^p$$
$$\leq 2^p \left(\frac{1}{n} \sum_{i=1}^n d_0(x_i, y_i)^p + \frac{1}{n} \sum_{i=1}^n d_0(x_{L(i)}, y_{L(i)})^p\right)$$
$$\leq 2^p (2\tau_D + 1)(d_1^{(p)}(\xi, \eta))^p,$$

using for the last inequality that any point of a point pattern in $(\mathbb{R}^D, |\cdot|)$ can be nearest neighbor to at most $\tau_D$ other points (see Zeger and Gersho (1994), Theorem 1). The factor $2^p$ is obviously unnecessary if $p = 1$. In Schuhmacher (2005a) a Lipschitz constant of $\tau_D + 1$ was obtained for $p = 1$ by a more complicated proof. □

## Acknowledgements

The author would like to thank the referees for many helpful suggestions that considerably improved the presentation of the paper. This work was supported by the Swiss National Science Foundation under Grant No. PBZH2-111668.